\documentclass{article}
\usepackage{amsmath}

\newtheorem{thm}{Theorem}[section]
\newtheorem{cor}[thm]{Corollary}
\newtheorem{prop}[thm]{Proposition}

\newtheorem{lem}[thm]{Lemma}

\newtheorem{ex}{Example}[section]

\newcommand{\be}{\begin{equation}}
\newcommand{\ee}{\end{equation}}
\newcommand{\ben}{\begin{enumerate}}
\newcommand{\een}{\end{enumerate}}
\newcommand{\pa}{{\partial}}

\newcommand{\qed}{\hspace*{\fill}Q.E.D.}  %Use at end of proof
\title{\Large On  Sprays with Vanishing $\chi$-Curvature}

\author{Zhongmin Shen\footnote{supported in part  by a NSFC grant (no: 11671352)}}

\date{January 1, 2020}

\begin{document}

\maketitle

\begin{abstract} Every Riemannian metric or Finsler metric on a manifold induces a spray via its geodesics.  
In this paper, we discuss  several expressions for the $\chi$-curvature of a spray. We show that the sprays obtained by a projective deformation using the S-curvature always have vanishing $\chi$-curvature. Then we establish the Beltrami Theorem for sprays with $\chi=0$.

{\bf MR(2000) subject classification: } 53C60, 53B40
\end{abstract}

\section{Introduction}

A spray $G$ on a manifold $M$ is a special vector field on the tangent bundle $TM$. 
In a standard local coordinate system $(x^i, y^i)$ in $TM$, a spray $G$ can be  expressed by 
\[  G = y^i \frac{\pa }{\pa x^i} - 2 G^i  \frac{\pa }{\pa y^i},\]
where $G^i = G^i(x,y) $ are local $C^{\infty}$ functions on non-zero vectors with the following homogeneity: $G^i(x, \lambda y)= \lambda^2 G^i(x,y),$ $ \forall \lambda >0.$
Every  Finsler metric induces a spray on a manifold. Some  geometric quantities of a Finsler metric are actually defined by the induced spray only. These quantities are extremely interesting to us.

For a spray $G$ on a manifold  $M$, with the Berwald connection,
 we define two key quantities: the Riemann curvature tensor $R^{\ i}_{j \ kl}$ and the Berwald curvature tensor $B^{\ i}_{j \ kl}$ (see \cite{Sh}).
Certain averaging process gives rise to various notions of  Ricci curvature tensor. One of them is the Ricci curvature tensor: 
$ {\rm Ric}_{jl} :=\frac{1}{2} \{  R^{\ m}_{j \ m l} + R^{\ m}_{l \ m j}\} $  (\cite{LiSh1}). 
The well known Ricci curvature 
${\rm Ric} := {\rm Ric}_{jl} y^jy^l = R^{\ m}_{j \ ml}y^j y^l$ has been studied for a long time by many people. 
Besides these quantities, we have another important quantity which is expressed in terms of the vertical derivatives of the Riemann curvature. It is the so-called $\chi$-curvature defined by
\be
  \chi_k:= -\frac{1}{6} \Big \{ 2 R^m_{\ k\cdot m} + R^m_{\ m \cdot k } \Big \}  .\label{chi_def}
\ee
where $R^i_{\ k} = y^j R^{\ i}_{j \ kl}y^l$. The $\chi$-curvature  can be expressed in several forms.
For an arbitrary volume form $dV$, 
\be
\chi_k = \frac{1}{2} \Big \{ S_{\cdot k|m} y^m-S_{|k} \Big \}, \label{chi_S}
\ee
where $S= S_{(G, dV)}$ is the S-curvature of $(G, dV)$  (\cite{Sh1}). 
For a spray induced by a Finsler metric,  the $\chi$-curvature can be expressed by
\be
  \chi_k = \frac{1}{2} \Big \{ I_{k|p|q} y^py^q + I_m R^m_{\ k}  \Big\}, \label{chi_I}
\ee 
where $I_k := g^{ij} C_{ijk}$ denotes the mean Cartan torsion  (\cite{Sh0} \cite{ChSh}). 
These are three typical expressions for the mysterious quantity $\chi$.
 In this paper, we shall focus on sprays with $\chi=0$.

For a spray $G$ on a manifold $M$, in  the  projectively equivalent class of $G$, there is always a spray  with $\chi=0$.  More precisely, for any volume form $dV$ on $M$, we may construct a spray $\hat{G}$  by a projective change:
\[  \hat{G}^i := G^i -\frac{S}{n+1} y^i,\]
where $S$ is the S-curvature of $(G, dV)$. This spray $\hat{G}$ is invariant under a projective change with $dV$ fixed. 
 This projective deformation is first introduced in \cite{Sh}. 
 We prove the following

\begin{thm}\label{thm1.1}
Let $G$ be a spray  on a manifold $M$. For any volume form $dV$, the  spray $\hat{G}$  associated with $(G, dV)$ has vanishing $\chi$-curvature,  $\hat{\chi}=0$.  
\end{thm}

Note that $\hat{G}$ is projectively equivalent to $G$. Hence if $G$ is of scalar curvature,
 then $\hat{G}$ is of scalar curvature too. Hence it is of isotropic curvature since  $\hat{\chi}=0$. Thus $\hat{G}$ must be  of isotropic curvature. 
We obtain the following

\begin{cor}\label{cor1.2}
Let $G$ be a spray  of scalar curvature on a manifold $M$. For any volume form $dV$, the  spray $\hat{G}$  associated with $(G, dV)$ must be of isotropic curvature. 
\end{cor}

The well-known Beltrami Theorem in Riemannian geometry says that for two projectively equivalent Riemannian metrics $g_1, g_2$,  the metric $g_1$ is of constant curvature if and only if $g_2$ is of constant curvature. In particular, if a Riemannian metric $g$ is locally projectively flat, then it is of constant curvature since $g$ is locally projectively equivalent to the standard Euclidean metric. This theorem can be extended to sprays with  $\chi=0$.

\begin{thm}\label{thm1.3} For two projectively equivalent sprays  $G_1, G_2$ with  $\chi=0$,   $G_1$ is of isotropic curvature  if and only if $G_2$ is of isotropic curvature. In particular, if a spray $G$ is locally projectively flat with $\chi=0$, then it is of isotropic curvature. 
\end{thm}

Sprays or Finsler metrics with $\chi=0$ deserve further study. Spherically symmetric metrics with $\chi=0$ have been studied in \cite{Zhu}.

\bigskip
\noindent
{\bf Acknowledgment}: The primary version of this note is part of my lectures during the summer school in 2018 in Xiamen University, China.

\section{Preliminaries}

A spray  $G$ on a manifold $M$ is a vector field on the tangent bundle $TM$  which is locally expressed in the following form
\[  G = y^i \frac{\pa }{\pa x^i} - 2 G^i  \frac{\pa }{\pa y^i},\]
where $G^i = G^i(x,y) $ are local $C^{\infty}$ function on $TU \equiv U \times R^n$,
\[ G^i(x, \lambda y ) = \lambda^2 G^i(x, y), \  \ \ \ \ \lambda >0.\]
Put
\[ N^i_j := \frac{\pa G^i}{\pa y^j}, \ \ \ \ \ \Gamma^i_{jk} = \frac{\pa^2 G^i}{\pa y^j \pa y^k}.\]
Let $\omega^i := dx^i$ and $\omega^{n+i} := dy^i + N^i_j dx^j$ and $\omega_j^{\ i} := \Gamma^i_{jk} dx^k$.
We have 
\[ d\omega^i = \omega^j \wedge \omega_j^{ \ i} .\]
Put 
\[  \Omega_j^{\ i} := d\omega_j^{\ i} -\omega_j^{\ k}  \wedge \omega_k^{\ i}.\]
We obtain two quantities $R$ and $B$:
\[
 \Omega_j^{\ i} =\frac{1}{2} R^{\ i}_{j \ kl} \omega^k \wedge \omega^l - B^{\ i}_{j \ kl} \omega^k \wedge \omega^{n+l},\]
where $R^{\ i}_{j \ kl} +  R^{\ i}_{j \ lk}=0$. 
\[  R^{\ i}_{j \ kl} =  {\delta \Gamma^i_{jl}\over \delta x^k} - {\delta \Gamma^i_{jk}\over \delta x^l}
+\Gamma^i_{ks}\Gamma^s_{jl} - \Gamma^s_{jk}\Gamma^i_{ls},  \]
\be
 B^{\ i}_{j \ kl} = \frac{\pa \Gamma^i_{kl}}{\pa y^j}.    \label{Bcurvature}
\ee

We have the first set of Bianchi identities
\begin{eqnarray}
&& R^{\ i}_{j \ kl} + R^{\ i}_{k \ lj} + R^{\ i}_{l\ jk} =0\\
&& B^{\ i}_{j \ kl} = B^{\ i}_{k \ jl}.
\end{eqnarray}
In fact $B^{\ i}_{j \ kl}$ is symmetric in $j, k, l$ and $y^j B^{\ i}_{j \ kl} =0$. 
Put  
\[ R^i_{\ kl} := y^jR^{\ i}_{j \ kl}, \ \ \ \ \  R^{\ i}_{j \ k}:= R^{\ i}_{j \ kl}y^l, \ \ \ \ \   R^i_{\ k} := y^j R^{\ i}_{j \ kl} y^l.\]
The two-index  Riemann curvature tensor $ R^i_{\ k} $ and the four-index Riemann curvature tensor $R^{\ i}_{j \ kl}$  determine each other by the following identity:
\be
R^{\ i}_{j \ kl} = \frac{1}{3} \Big \{ R^i_{\ k\cdot l\cdot j} - R^i_{\ l\cdot k\cdot j} \Big \},  \label{RRR}
\ee
We also have 
\begin{eqnarray}
R^{\ i}_{j \ k} & = & \frac{1}{3} \Big \{ 2 R^i_{\ k\cdot j}+R^i_{\ j \cdot k} \Big \}, \label{RRR2}\\
R^i_{\ kl} & = & \frac{1}{3} \Big \{ R^i_{\ k\cdot l} - R^i_{\ l\cdot k} \Big \},\label{RRR3}
\end{eqnarray}
where $T^{*}_{\ *\cdot k} $ is the vertical covariant derivative, namely, $T^*_{*\cdot k} = \frac{\pa}{\pa y^k} ( T^*_{\ *})$.

Further covariant derivatives yield the second set of Bianchi identities: 
\begin{eqnarray}
&&R^{\ i}_{j \ kl|m}+ R^{\ i}_{j \ lm|k} + R^{\ i}_{j \ mk|l}\nonumber\\
&& \hspace{1 cm}   + B^{\ i}_{j \ mp}R^p_{ \ kl} + B^{\ i}_{j \ lp}R^p_{\ mk} + B^{\ i}_{j \ kp}R^p_{\ lm} =0\label{Bi1}\\
&& R^{\ i}_{j \ kl\cdot m} = B^{\ i}_{j\ ml|k} - B^{\ i}_{j\ km|l}\label{Bi2} \\
&& B^{\ i}_{j \ kl\cdot m} = B^{\ i}_{j \ km\cdot l}. \label{Bi3}
\end{eqnarray}
Contracting (\ref{Bi1}) with $y^j$ yields
\be
R^i_{\ kl|m} + R^i_{\ lm|k}+ R^i_{\ mk|l} =0. \label{Bi4}
\ee
Contracting (\ref{Bi4}) with $y^l$ yields
\be
R^i_{\ k|m} - R^i_{\ m|k} + R^i_{\ mk|l}y^l =0. \label{RRRR}
\ee

\section{The $\chi$-curvature}

The   $\chi$-curvature   defined by the Riemann curvature tensor  in (\ref{chi_def}) can be expressed in several ways.

\begin{lem}\label{lem3.1}
\be
\chi_k =  -\frac{1}{2} R^{\ m}_{m \ k} = - \frac{1}{2} R^{\ m}_{m \ kl}y^l.  \label{XR}
\ee
\end{lem}
{\it Proof}:  It follows from (\ref{RRR2}). 
\qed

\bigskip
Lemma \ref{lem3.1} tells us that if $R^{\ m}_{m\  k}=0$, then $\chi=0$.

Put
\be
T^i_{\ k} : = R^i_{ \ k} - \Big\{ R \delta^i_{\ k} - \frac{1}{2} R_{\cdot k} y^i \Big \},  \label{Tcurvature}
\ee
where $R:= \frac{1}{n-1} R^m_{\ m}$. 
By definition, $G$ is of isotropic curvature if $T^i_{\ k}=0$.  
Note that 
\[  {\rm trace} (T) := T^m_{\ m} =0.\]

By a direct computation, we can obtain another expression for $\chi_k$.
\begin{lem}\label{lem3.2}
\be
 \chi_k = -\frac{1}{3} T^m_{\ k\cdot m}.
\ee
 \end{lem}

\bigskip
Lemma \ref{lem3.2} tells us that if $G$ is of isotropic curvature, then  $\chi=0$. 

Recall the definition of the Weyl curvature
\be
W^i_{\ k} := A^i_{\ k} - \frac{1}{n+1} A^m_{\ k \cdot m} y^i,
\ee 
where $A^i_{ \ k} := R^i_{\ k} - R \delta^i_{\ k}$.  Clearly,
\[ W^m_{\ k \cdot m}=0.\]
We obtain a nice formula for the Weyl curvature. 

\begin{lem}
The Weyl curvature is given by
\be
W^i_{\ k} = R^i_{\ k} -\Big \{ R\delta^i_{\ k} - \frac{1}{2} R_{\cdot k} y^i \Big \} +\frac{3}{n+1} \chi_k  y^i.  \label{WRX}
\ee
\end{lem}
{\it Proof}: One can easily rewrite $W^i_{\ k}$ as 
\[ W^i_{\ k} = R^i_{\ k} - \Big \{ R\delta^i_{\ k} - \frac{1}{2}R_{\cdot k} y^i \Big \} -\frac{1}{2(n+1)}
\Big \{ 2R^m_{\ k \cdot m} +(n-1) R_{\cdot k}\Big \} y^i.\]
By (\ref{XR}), we prove the lemma.
\qed

\bigskip

Given a volume $dV= \sigma (x) dx^1 \cdots dx^n$, the S-curvature of $(G, dV)$ is defined by 
\[   S := \Pi- y^m \frac{\partial }{\partial x^m} \Big ( \ln \sigma \Big ).\]
We have the following expression for $\chi$. 
\begin{lem}(\cite{LiSh1})
\be 
 \chi_k  = \frac{1}{2} \Big \{ S_{\cdot k |m} y^m - S_{|k} \Big \}.  \label{XS}
\ee
\end{lem}

In local coordinates, by (\ref{XS}), one can easily get  
\be
\chi_k = \frac{1}{2} \Big \{  \Pi_{x^m y^k} y^m - \Pi_{x^k} -2 \Pi_{y^k y^m} G^m \Big \},  \label{Slocal}
\ee
where $ \Pi := \frac{\partial G^m}{\partial y^m} $.   Clearly, $\chi$ is independent of $dV$.

\bigskip

\section{Sprays with $\chi=0$}

A spray is said to be  {\it  $S$-closed} if in local coordinates, $\Pi = \frac{\pa G^m}{\pa y^m} $ is a closed local $1$-form. 
The spray induced by a Riemannian metric $g=g_{ij}(x)y^iy^j$ is S-closed. 
In fact 
\be
 \Pi = y^k \frac{\pa }{\pa x^k} \Big [ \ln  \sqrt{\det(g_{ij}(x) ) } \Big ] .\label{Piclosed}
\ee
 By (\ref{Piclosed}), for any volume form $dV=\sigma(x) dx^1\wedge \cdots \wedge dx^n$,  the S-curvature of $(G, dV)$ is a closed $1$-form,
\[  S = y^k \frac{\pa  }{\pa x^k}[\ln \varphi(x)] ,\]
where $\varphi ( x) =   \sqrt{\det(g_{ij}(x) ) }/\sigma(x)$.

We have the following

\begin{prop}\label{propS=closed} If a spray is S-closed, then $\chi=0$.
In particular, if for some volume form $dV = \sigma dx^1 \cdots dx^n$, the S-curvature of $(G, dV)$ is a closed $1$-form, then $\chi=0$.
\end{prop}
{\it Proof}: 
By assumption,
\[ S = \Pi - y^m \frac{\pa }{\pa x^m} (\ln \sigma)  = \eta _k y^k, \]
with $(\eta_k)_{x^l} = (\eta_l)_{x^k}$.  Then by (\ref{Slocal}), $\chi_k=0$.
\qed

\bigskip

Let $\tilde{F}$  be a Finsler metric  and $G$ be a spray on a manifold $M$. 
The spray coefficients $\tilde{G}^i$ of $\tilde{F}$ can be expressed  as follows
\be
\tilde{G}^i = G^i + P y^i 
+ \frac{1}{2} \tilde{F} \tilde{g}^{ik} \Big \{ \tilde{F}_{\cdot k |m} y^m -\tilde{F}_{|k} \Big \}.
\ee 
where $P= \tilde{F}_{|m}y^m/(2\tilde{F})$. 
Thus  $\tilde{F}$ is projectively equivalent to $G$ if and only if 
\be
\tilde{F}_{\cdot k |m} y^m -\tilde{F}_{|k} =0. 
\ee
This is the generalized version of the famous Rapcs\'{a}k Theorem. 
By (\ref{XS}), we obtain the following

\begin{thm}\label{thm4.3} Let $G$  be a spray  with $\chi=0$ and $ dV$  be a volume form. If for  the S-curvature  $S$ of $(G, dV)$,   $\tilde{F}=|S|$ is a Finsler metric, then it is projectively equivalent to $G$.  
\end{thm}

\section{Sprays of isotropic curvature}

A spray $G$ is said to  be {\it of scalar curvature} if 
\be
R^i_{ \ k} = R \delta^i_{\ k} - \tau_k y^i , \label{RRt}
\ee
where $\tau_k $ is a positively homogeneous function of degree one with $\tau_k y^k =  R$. 
This is equivalent to that $W^i_{\ k} =0$. By  (\ref{WRX}), we see that (\ref{RRt}) is equivalent to
the following
\be
R^i_{ \ k} = R \delta^i_{\ k} -\frac{1}{2} R_{\cdot k} y^i -\frac{3}{n+1}\chi_k y^i. \label{RRCS}
\ee

\bigskip

The $\chi$-curvature characterizes sprays of isotropic curvature among  sprays of scalar curvature.  By (\ref{RRCS}), we obtain the following

\begin{thm}\label{propsi} (\cite{LiSh2}) \label{lemchi=0} Let $G$ be a spray of scalar curvature. $G$ is of isotropic curvature if and only if $\chi=0$.
\end{thm}

\noindent
{\it Proof of Theorem \ref{thm1.3}}: If $G_1$ is of isotropic curvature, then $G_2$ is of scalar curvature by the projective equivalence. Since $\chi=0$, we see that $G_2$ is of isotropic curvature by Proposition \ref{propsi}.
\qed

\bigskip
If $G$ is of isotropic curvature, then 
\[
R^{\ i}_{j \ kl} = \frac{1}{2} \Big \{ R_{\cdot l \cdot j} \delta^i_{\ k} - R_{\cdot k \cdot j} \delta^i_{\ l} \Big \}.
\]
\[
R^i_{\ kl} =\frac{1}{2} \Big \{  R_{\cdot l} \delta^i_{\ k} -R_{\cdot k}\delta^i_{\ l} \Big \}.
\]

Assume that $G$ is of isotropic curvature.
By (\ref{Bi1}), we obtain
\be
(R_{\cdot l\cdot j|m}-R_{\cdot m \cdot j |l} )\delta^i_{\ k} 
+ (R_{\cdot m\cdot j|k} - R_{\cdot k\cdot j|m} )\delta^i_{\ l} 
+ (R_{\cdot k \cdot j |l} -R_{\cdot l \cdot j |k} )\delta^i_{\ m} =0.
\ee
This yields
\be
(R_{\cdot l|m}-R_{\cdot m|l})\delta^i_{\ k} 
+(R_{\cdot m|k} - R_{\cdot k|m} )\delta^i_{\ l} 
+ (R_{\cdot k|l} - R_{\cdot l|k} )\delta^i_{\ m}=0.  \label{Bi5}
\ee
Contracting  (\ref{Bi5}) with $y^m$ yields
\be
 (R_{\cdot l|m}y^m- 2 R_{|l})\delta^i_{\ k} + (2 R_{|k}-R_{\cdot k|m}y^m) \delta^i_{| k}  
+ (R_{\cdot k|l}-R_{\cdot l|k}) y^i =0. \label{RR1}
\ee
Taking trace $ i=k$  in (\ref{RR1}), we obtain
\be
 (n-2) (R_{\cdot l|m}y^m - 2 R_{|l} ) =0.\label{RR2}
\ee

\begin{thm} \label{thmdual}
If $G$ is  an $n$-dimensional spray of isotropic curvature $R$  ($n\geq 3$), then $R$ satisfies
\be
 \frac{1}{2} R_{\cdot  l |m} y^m - R_{|l}=0.  \label{RRik}\ee
\end{thm}
{\it Proof}: By assumption $n\geq 3$,  we obtain from (\ref{RR2}), 
\[ R_{|l} - \frac{1}{2} R_{\cdot l|m}y^m =0.\]
\qed

\bigskip

For a spray $G$, we introduce a new quantity  $\eta = \eta_k  d^k$, 
\be
 \eta_k:=  \frac{1}{2} R_{\cdot k|m}y^m - R_{|k},  \label{eta}
\ee
where $R:=\frac{1}{n-1} {\rm Ric} $.

For a spray of isotropic curvature $R$ on $n$-dimensional manifold $M$ ($n\geq 3$),  
By Theorem \ref{thmdual}, $\eta =0$.

\bigskip
Let $L:=\tilde{F}^2$ be a Finsler metric and $G$ a spray on a manifold.   The spray coefficients of $\tilde{F}$ can be expressed as 
\be
\tilde{G}^i=G^i + \frac{1}{4} 
\tilde{g}^{ik} L_{|k}  + \tilde{g}^{ik} \Big \{ \frac{1}{2} L_{\cdot k | m} y^m - L_{|k} \Big \}.  \label{Gdual}
\ee
$L:=\tilde{F}^2$  is said to be {\it dually equivalent to} $G$ if 
\be
\tilde{G}^i=G^i + \frac{1}{4} 
\tilde{g}^{ik} L_{|k}.
\ee
This is  equivalent to 
\be
 \frac{1}{2} L_{\cdot k | m} y^m - L_{|k}  =0.  \label{LLL}
\ee

\bigskip

For a spray $G$ on an $n$-dimensional manifold $M$ with isotropic scalar curvature $R$.  Assume that $R$ is a Finsler metric,  by Theorem \ref{thmdual}, one can see that 
$R$ is dually equivalent to $G$.

\section{Projective change by the S-curvature}
Let $G$ be a spray  and $dV$ be a volume form on an $n$-dimensional manifold $M$.  We deform $G$ to another spray $\hat{G}$ by 
\[  \hat{G}^i :=  G^i  - \frac{S}{n+1} y^i,\]
where $S $ denotes the S-curvature of $(G, dV)$.  From  the definition, we see that $\hat{G}$ is projectively equivalent to $G$. 

\bigskip

\begin{lem}  Let $G$ be a spray and $dV$ a volume form on a manifold $M$.  Let $\hat{G}$ be  the spray associated with $(G, dV)$. Then  the S-curvature  of $(\hat{G}, dV)$ vanishes. Hence, 
$\hat{\chi} =0$. 
\end{lem}
{\it Proof}: Recall
\[ \hat{\chi}_k = \frac{1}{2} \Big \{ \hat{S}_{|m \cdot k} y^m - \hat{S}_{|k} \Big \}.\]
On the other hand,  $\hat{G}^i = G^i + Py^i$ with $P = - \frac{S}{n+1}$. Thus 
\[ \hat{S}= S +(n+1) P = 0.\]
This yields that $ \hat{\chi}=0$. 
\qed

\begin{lem}
If $G_1$ and $G_2$ are two projectively equivalent sprays on a manifold $M$, then for any volume form $dV$, the spray $\hat{G}_1$  associated with 
$(G_1, dV)$  and  $\hat{G}_2$ associated with $(G_2, dV)$ are equal, i.e., $\hat{G}_1 = \hat{G}_2$. 
\end{lem}
{\it Proof}: It is easy to see that  if $ G^i_1 = G^i_2 + P y^i$, then 
\[  S_1 = S_2 + (n+1) P.\]
Then 
\begin{eqnarray*}
\hat{G}_1^i & = &  G^i_1 - \frac{S_1}{n+1} y^i \\
& = &  [G^i_2 + Py^i] - \frac{S_2+(n+1) P}{n+1} y^i \\
& = &  G^i_2 - \frac{S_2}{n+1} y^i = \hat{G}_2.
\end{eqnarray*}
\qed

\bigskip

\noindent
{\it Proof of Corollary \ref{cor1.2}}: First by definition, $\hat{G}$ is projectively equivalent to $G$. Thus $\hat{G}$ is of scalar curvature.  Since $\hat{\chi}=0$, by Lemma \ref{lemchi=0},  we see that $\hat{G}$ is of isotropic curvature. \qed

\bigskip

By the above lemma, any geometric quantity of $\hat{G}$ is a projective invariant of $G$ with respect to a fixed volume form $dV$.  Further, if the geometric quantity of $\hat{G}$ is independent of  the volume form $dV$, then the quantity is a projective quantity of $G$.

\begin{lem}  Let $G$ be a spray and $dV$  a volume form on a manifold $M$.  For the spray $\hat{G}$ associated with $(G, dV)$, the Riemann curvature of $\hat{G} $ is given by
\be
\hat{R}^i_{\ k} =  R^i_{\ k} + \tau \delta^i_{\ k} - \frac{1}{2} \tau_{\cdot k} y^i +\frac{3 \chi_k}{n+1} y^i, \label{RRtaux}
\ee
where 
\be
\tau : = \Big ( \frac{S}{n+1} \Big )^2 + \frac{1}{n+1} S_{|m}y^m.
\ee
\end{lem}
{\it Proof}: By a direct argument. 
\qed

\bigskip

 By (\ref{RRtaux}), we get the projective Ricci curvature tensor $\widehat{\rm Ric}_{jl}
:=\frac{1}{2} \{  \hat{R}^{\ m}_{j \ ml} +\hat{R}^{\ m}_{l \ mj} \}$ and the projective Ricci curvature $\widehat{\rm Ric} := \widehat{\rm Ric}_{jl}y^jy^l$.  

\begin{eqnarray}
\widehat{\rm Ric}_{jl} & = & {\rm Ric}_{jl} +\frac{n-1}{2} \tau_{\cdot j\cdot l} +H_{jl}, \label{hatRicij}\\
\widehat{\rm Ric} & =  & {\rm Ric} + (n-1) \tau,\label{hatRic}
\end{eqnarray}
where  $\widehat{\rm Ric} = \widehat{\rm Ric}_{jl} y^j y^l$ is the Ricci curvature of $\hat{G}$ and 
\[
 H_{ij} := \frac{1}{2}\Big \{ \chi_{i\cdot j}+\chi_{j\cdot j} \Big \}.
\]

\bigskip
It is natural to consider other quantities of $\hat{G}$, such as the Berwald curvature defined in (\ref{Bcurvature})
and the T-curvature defined in (\ref{Tcurvature})
\[  \hat{B}^{\ i}_{j \ kl} = \frac{\pa^3 \hat{G}^i}{\pa y^j \pa y^k \pa y^l}.\]
\[   \hat{T}^i_{\ k} = \hat{R}^i_{\ k} - \Big \{ \hat{R}\delta^i_{\ k} - \frac{1}{2} \hat{R}_{\cdot k} y^i \Big \}.\]
Clearly, $\hat{B}$ and $\hat{T}$ are projective invariants with a fixed volume form $dV$.  We have the following

\begin{prop}  Let $G$ be a spray on  a manifold  and $\hat{G}$  a spray  associated with $(G, dV)$ for some volume form $dV$. Then the Berwald curvature $\hat{B}$ and $\hat{T}$ are independent of $dV$, hence they are projective invariants of $G$. In fact 
$ \hat{B} = D$ is the Douglas curvature and $\hat{T} = W$ is the Weyl curvature  of $G$. 
\end{prop}

Here we provide another description of the Douglas curvature and the Weyl curvature of a spray.

Let $G$ be a spray and $\hat{G}$ be the spray associated with $(G, dV)$ for some volume form $dV$. 
Let $\hat{\eta}$ be the quantity of $\hat{G}$  defined in (\ref{eta}).  Then $\hat{\eta}$ is a projective invariant of $G$ for a fixed  volume form $dV$.  In fact, $\hat{\eta}  = {\bf W}^o$ the so-called  {\it Berwald-Weyl curvature} (\cite{Sh}).  If  $G$ is of scalar curvature, then
$\hat{G}$ is of isotropic curvature. Thus  $\hat{\eta}=0$ when $n=\dim M \geq 3$  by Theorem \ref{thmdual}.

\begin{prop}Let $G$ be a spray on  a manifold  and $\hat{G}$  a spray  associated with $(G, dV)$ for some volume form $dV$.  Assume that $G$ is of scalar curvature. 
Then the projective invariant  $\hat{\eta} =0$ in dimension $n\geq 3$.
\end{prop}

\section{Examples}
In this section, we shall give some sprays of isotropic curvature.

\bigskip
\begin{ex}{\rm
Let 
$ F =\alpha+\beta$ be a Randers metric on an $n$-dimensional manifold $M$, 
where $\alpha =\sqrt{a_{ij}(x)y^iy^j}$ is a Riemannian metric and $\beta = b_i (x) y^i$ is a $1$-form on $M$. Let $\nabla \beta =b_{i|j} y^i dx^j $ denote the covariant derivative of $\beta$ with respect to $\alpha$. Let 
\[ r_{ij}:= \frac{1}{2} (b_{i|j} + b_{j|i} ), \ \ \ \  s_{ij}:= \frac{1}{2} ( b_{i|j}-b_{j|i} ), \ \ \ \ s_j := b^i s_{ij},\]
\[ q_{ij}:= r_{im}s^m_{\ j}, \ \ \ t_{ij}:=s_{im}s^m_{\ j}, \ \ \ \ t_j:=b^i t_{ij}.\]

Let 
\be
 \hat{G}^i: = G^i_{\alpha} + \alpha s^i_{\ 0}.\label{RandersG}
\ee
In fact $\hat{G}$ is the spray  associated with $(G, dV_{\alpha})$. 
It is proved that 
$\hat{G}$ is of scalar curvature if and only if the Riemann curvature $\bar{R}^i_{\ k}$ of $\alpha$ and the covariant derivatives of  $\beta$ satisfy the following equations  (\cite{ShYi})
\begin{eqnarray}
\bar{R}^i_{\ k} & = & \kappa  \Big \{ \alpha^2 \delta^i_k -  y_k y^i \Big \}\nonumber\\
&& +\alpha^2 t^i_{\ k}+ t_{00} \delta^i_k-t_{k0}y^i
  -  t^i_{\ 0} y_k - 3s^i_{\ 0} s_{k0} , \label{eqW1AAA}\\
s_{ij|k} &  =  &  \frac{1}{n-1} \Big \{a_{ik} s^m_{\ j|m} - a_{jk} s^m_{\ i|m}  \Big \}. \label{eqW2**AA}
\end{eqnarray}
where $\kappa =\kappa (x)$ is a scalar function on $M$.  
In this case,  $\hat{G}$ is actually of isotropic curvature. $\hat{R}^i_{\ k} = \hat{R}\delta^i_k -\frac{1}{2} \hat{R}_{\cdot k} y^i$. 
By a simple computation, we obtain a formula for $\hat{R}:= \frac{1}{n-1} \widehat{\rm Ric}$:
\[  \hat{R} = \kappa \alpha^2 + t_{00}+ \frac{2}{n-1} \alpha s^m_{\ 0|m}.\]

}
\end{ex}

. 

\begin{ex}{\rm 
Consider a spray on an open subset $U\subset R^2$,
\[ G = y^1 \frac{\pa }{\pa x^1}+ y^2 \frac{\pa }{\pa x^2} - 2 G^1 \frac{\pa }{\pa y^1}-2 G^2 \frac{\pa }{\pa y^2},\]
where 
\begin{eqnarray*}
G^1 & = &  B (y^1)^2 + 2 C y^1y^2 + D (y^2)^2 +\frac{1}{3} ( f_{x^1} (y^1)^2 +f_{x^2} y^1y^2)\\
G^2 & = & - A (y^1)^2 -2 B  y^1y^2 - C (y^2)^2 +\frac{1}{3} (f_{x^1} y^1y^2 +f_{x^2} (y^2)^2).
\end{eqnarray*}
where 
\[  A = A(x^1,x^2), \ \  B = B(x^1,x^2), \ \ C = C(x^1,x^2), \ \ D = D(x^1, x^2), \ \ f = f(x^1, x^2)\]
are $C^{\infty}$ functions on $U$.  The geodesics are the graphs of $x^2 = \phi (x^1) $
\[  \phi''  = 2 A(x^1, \phi) + 6 B(x^1, \phi) \phi'+ 6 C(x^1, \phi) (\phi')^2  + 2 D(x^1, \phi) (\phi')^3.\]

 We have
\[ \Pi =\frac{\pa G^m}{\pa y^m} = f_{x^1} y^1 + f_{x^2} y^2.\]
Thus $\chi_k =0$.
Further computation shows that $G$ is of isotropic curvature. }

\end{ex}

\bibliographystyle{plain}

\begin{thebibliography}{lbl}



\bibitem{ChSh}X. Cheng and Z. Shen, {\it  Finsler Geometry --- An approach via Randers spaces}, Springer-Verlag, (2012)

\bibitem{LiSh1} B. Li and Z. Shen, {\it  Ricci curvature tensor and non-Riemannian quantities}, Canadian Mathematical Bulletin, 58(2015), 530-537.

\bibitem{LiSh2} B. Li and Z. Shen, {\it  On sprays of isotropic curvature}, International Journal of Mathematics, 29 (2018), https://doi.org/10.1142/S0129167X18500039

\bibitem{Sh0}  Z. Shen, {\it Finsler manifolds with nonpositive flag curvature and constant S-curvature}, Mathematische Zeitschrift, {\bf 249}(2005), 625-639.

\bibitem{Sh1}
 Z. Shen, {\it On some non-Riemannian quantities in Finsler geometry}, Canad. Math. Bull.  {\bf 56}(2013), 184-193. 

\bibitem{Sh} Z. Shen, {\it Differential Geometry of Spray and Finsler Spaces},
Kluwer Academic Publishers, 2001.

\bibitem{ShYi} Z.Shen and G. C. Yildirim, {\it A characterization of Randers metrics of scalar flag curvature}, Recent Developments in Geometry and Analysis, Advanced Lectures in Mathematics {\bf 23} (2013), 345-358. 


\bibitem{Zhu}  H. Zhu, {\it On a class of Finsler metrics with special curvature properties}, 
Balkan Journal of Geometry and Its Applications, {\bf 23} (2018), 97-108.

\end{thebibliography}

\vspace{0.6cm}

\noindent Zhongmin Shen

\noindent Department of Mathematical Sciences, Indiana
University-Purdue University Indianapolis, IN 46202-3216, USA.

\noindent \verb"zshen@math.iupui.edu"

\end{document}